\documentclass[12pt]{amsart}

\usepackage{setspace}   %Requested by Akhil
\usepackage[usenames]{color}
\usepackage{fullpage,url,mathrsfs,stmaryrd}
\usepackage{amsmath,amsthm,amssymb,mathrsfs,stmaryrd,color}
\usepackage[utf8]{inputenc}
\usepackage[T1]{fontenc}
\usepackage{eucal}

\usepackage{relsize}
\usepackage[bbgreekl]{mathbbol}
\usepackage{amsfonts}

\DeclareSymbolFontAlphabet{\mathbb}{AMSb}
\DeclareSymbolFontAlphabet{\mathbbl}{bbold}

\newcommand{\BG}{{\mathbb{G}}}

\usepackage{amssymb}
\usepackage[all]{xy}
\usepackage{mathrsfs}
\usepackage{enumerate}
\usepackage{amscd}

\usepackage{bbm}

\DeclareMathOperator{\KS}{{KS}}
\DeclareMathOperator{\pCurv}{{p-Curv}}

\newcommand{\EEnd}{\underline{\on{End}}}

\newcommand{\HHom}{\underline{\on{Hom}}}

\DeclareMathOperator{\rank}{{rank}}

\newcommand{\cA}{{\mathcal A}}
\newcommand{\cB}{{\mathcal B}}

\newcommand{\cO}{{\mathcal O}}

\newcommand{\nc}{\newcommand}

\nc\wh{\widehat}

\nc\on{\operatorname}

\nc\Gr{\on{Gr}}

\nc\Fl{\on{Fl}}

\newtheorem{problem}[subsubsection]{Problem}

\theoremstyle{remark}

\newcommand{\BC}{{\mathbb{C}}}
\newcommand{\BF}{{\mathbb{F}}}
\newcommand{\BN}{{\mathbb{N}}}
\newcommand{\BQ}{{\mathbb{Q}}}
\newcommand{\BR}{{\mathbb{R}}}

\newcommand{\BZ}{{\mathbb{Z}}}

\DeclareMathOperator{\BT}{{BT}}

\newcommand{\limto}{{\displaystyle\lim_{\longrightarrow}}}
\newcommand{\rightlim}{\mathop{\limto}}

\newcommand{\leftlim}{\mathop{\displaystyle\lim_{\longleftarrow}}}
\newcommand{\limfromn}{\leftlim\limits_{\raise3pt\hbox{$n$}}}
\newcommand{\limton}{\rightlim\limits_{\raise3pt\hbox{$n$}}}

\newcommand{\rightlimit}[1]{\mathop{\lim\limits_{\longrightarrow}}\limits%
                    _{\raise3pt\hbox{$\scriptstyle #1$}}}

\newcommand{\leftlimit}[1]{\mathop{\lim\limits_{\longleftarrow}}\limits%
                    _{\raise3pt\hbox{$\scriptstyle #1$}}}

\newcommand{\iso}{\buildrel{\sim}\over{\longrightarrow}}

\DeclareMathOperator{\dR}{{dR}}

\DeclareMathOperator{\Ker}{{Ker}} 
\DeclareMathOperator{\im}{{Im}}

\DeclareMathOperator{\Spec}{{Spec}}

\theoremstyle{definition}

\numberwithin{equation}{section}

\newcommand{\Fr}{\operatorname{Fr}}

\begin{document}
\title[Toward Shimurian analogs of Barsotti-Tate groups]{Toward Shimurian analogs of Barsotti-Tate groups}
%\title[Shimurian generalizations of $\BT_n$]{Towards Shimurian generalizations of the stacks $\BT_n$}
%\title[On the stacks $\BT_1^G\otimes\BF_p$]{On the stacks $\BT_1^G\otimes\BF_p$}
\author{Vladimir Drinfeld}
\address{University of Chicago, Department of Mathematics, Chicago, IL 60637}
%\email{drinfeld@math.uchicago.edu}

%\begin{abstract}
%Let $G$ be a smooth group scheme over $\BF_p$ equipped with a $\BG_m$-action such that all weights of $\BG_m$ on $\Lie (G)$ are $\le 1$.
 %Let $\Disp_n^G$ be Eike Lau's stack of $n$-truncated $G$-displays (this is an algebraic $\BF_p$-stack).  
% In the case $n=1$ we introduce an algebraic stack equipped with a morphism to $\Disp_1^G$. We conjecture that if $G=GL(d)$ then the new stack is  canonically isomorphic to the reduction modulo $p$ of the stack of $1$-truncated Barsotti-Tate groups of height $d$ and  dimension $d'$, where $d'$ depends on the action of $\BG_m$ on $GL(d)$.
 
%We also discuss how to define an analog of the new stack for $n>1$ and how to replace $\BF_p$ by $\BZ/p^n\BZ$.
%\end{abstract}

%\keywords{Barsotti-Tate group, Shimura variety, display, $F$-zip, connection, $p$-curvature, syntomification}
%\subjclass[2010]{14F30}

\maketitle

These are the notes of the lecture that I gave on November 13, 2023 
at the Chinese University of Hong Kong. (During their preparation I was partially supported by NSF grant DMS-2001425.)

\medskip

The word ``toward'' from the title is important: I will be discussing an emerging theory (due to many people) rather than a theory which already exists. The word ``Shimurian'' refers to the notion of \emph{Shimura variety} (to be discussed in \S\ref{ss:Shimura} below). 
%The abbreviation ``BT'' stands for ``Barsotti-Tate''.

Fix a prime $p$ and a number $n\in\BN$. 

%\tableofcontents

\section{The notion of $\BT_n$ group}

$\BT_n$ is an abbreviation  for ``$n$-truncated Barsotti-Tate''.

\subsection{Grothendieck's definition of $\BT_n$ group}
Let $R$ be a (commutative) ring.

\subsubsection{Format of the definition}   \label{sss:format of definition}
A $\BT_n$ group over $R$ is a commutative finite group scheme $G$ over $R$ killed by $p^n$ with certain extra properties. The properties will be formulated 
in~\S\ref{sss:extra properties}. But first, let me explain what they are in a simple example.

\subsubsection{Example}  \label{ex:1}
Let $R$ be an algebraically closed field of characteristic $0$. Then a finite group scheme is the same as a finite group $G$. In this situation there is only one extra property: $G$~should be flat (or equivalently, free) as a $(\BZ/p^n)$-module. So $G\simeq (\BZ/p^n)^h$ for some number~$h$, which is called the \emph{height} of $G$. Note that
\begin{equation}  \label{e:h}
|G|=p^{nh},
\end{equation}
where $|G|$ is the order of $G$.

\subsubsection{Remark}
We will mostly be interested in the situation where $p$ is nilpotent in $R$. This situation is opposite to that of Example~\ref{ex:1}.

\subsubsection{The height of a $\BT_n$ group}
For any $R$, the height $h$ of a $\BT_n$ group $G$ over $R$ is defined by formula~\eqref{e:h}. (Both $|G|$ and $h$ are locally constant functions $\Spec R\to\BZ$.)

\subsubsection{The extra properties}   \label{sss:extra properties}
The reader may prefer to skip this subsection, in which we formulate the properties promised in \S\ref{sss:format of definition}.

(i) First, assume that $R$ is a field. Then commutative finite group schemes over $R$ form an abelian category, and our $G$ is a $(\BZ/p^n)$-module in this category. If $n>1$ then this module is required to be flat over $\BZ/p^n$. If $n=1$ this is automatic, but in the case that $R$ has characteristic $p$ and $n=1$ there is another requirement: namely, one requires $\Ker F=\im V$, where $F:G\to\Fr^*G$ is the Frobenius and $V:\Fr^*G\to G$ is the Verschiebung.

(ii) Now let $R$ be any ring. Then there are two requirements. First, for any field $R'$ and any homomorphism $f:R\to R'$, the base change of $G$ to $R'$ is required to be a $\BT_n$ group. Second, $G$ is required to be locally free; in other words, the coordinate ring of $G$ is required to be a projective $R$-module.

\section{Examples of $\BT_n$ groups}  \label{s:Examples}

\subsection{Two simple examples}  \label{ss:simple xamples}
\subsubsection{} The group scheme $\BZ/p^n$ is a $\BT_n$ group of height 1.

\subsubsection{}  Let 
\[
\mu_{p^n}:=\Ker (\BG_m\overset{p^n}\longrightarrow\BG_m),
\]
 where $\BG_m$ is the multiplicative group (viewed as a group scheme over $R$); in other words, 
$\mu_{p^n}$ is the group scheme over $R$ formed by $p^n$-th roots of $1$. Then $\mu_{p^n}$ is a $\BT_n$ group. Its order equals $p^n$, so its height equals $1$.

Let us look more attentively at $\mu_{p^n}$ in the case that $R$ is a field of characteristic $p$. The subscheme $\mu_{p^n}$ is defined by the equation $x^{p^n}-1=0$. In characteristic $p$ one has 
$$x^{p^n}-1=(x-1)^{p^n},$$
so $\mu_{p^n}$ has a single $R$-point (namely, $x=1$). Nevertheless, the order of $\mu_{p^n}$ equals $p^n$, as stated above. Indeed, the order of a finite scheme is defined to be the dimension of its coordinate ring, and the coordinate ring of $\mu_{p^n}$ is $R[x]/(x^{p^n}-1)$, so it has dimension $p^n$.
 
\subsection{The motivating (and typical) example}  \label{ss:motivating example}
\subsubsection{The example}  \label{sss:abelian schemes}
Let $A$ be an abelian scheme over $R$ (i.e., a family of abelian varieties parametrized by $\Spec R$). Then the group scheme
\begin{equation}  \label{e:abelian schemes}
A[p^n]:=\Ker (A\overset{p^n}\longrightarrow A)
\end{equation}
is a $\BT_n$ group of height $2d$, where $d$ is the dimension of $A$ over $R$. This explains why $\BT_n$ groups are interesting and useful. 

$\BT_n$ groups of the form \eqref{e:abelian schemes} are ``typical'' in the following sense: it is known that locally for the flat topology of $\Spec R$,
every $\BT_n$ group can be represented as a direct summand of a $\BT_n$ group of the form \eqref{e:abelian schemes}.

\subsubsection{The case $R=\BC$}
A $d$-dimensional abelian variety over $\BC$ is a complex torus of complex dimension $d$ and real dimension $2d$. So the group \eqref{e:abelian schemes} is isomorphic to $(\BZ/p^n)^{2d}$. Therefore it is a $\BT_n$ group of height $2d$, as already said in \S\ref{sss:abelian schemes}.

\subsubsection{The case $R=\bar\BF_p$, $d=1$}  \label{sss:elliptic curves}
Abelian varieties of dimension $1$ are called elliptic curves. If $A$ is an elliptic curve then the $\BT_n$ group $A[p^n]$ defined by \eqref{e:abelian schemes} has height $2$. It is known that $A[p^n]\simeq (\BZ/p^n)\oplus\mu_{p^n}$ for all but finitely many elliptic curves over $\bar\BF_p$. The exceptional elliptic curves are called supersingular. We will not discuss the structure of $A[p^n]$ for a supersingular elliptic curve $A$; let me just say that in this case the group of $\bar\BF_p$-points of $A[p^n]$ is zero.

\subsubsection{Remark}
If $G$ is a $\BT_n$ group over a ring $R$ then the fibers of $G$ over different points of $\Spec R$ can be quite different from each other (even if $\Spec R$ is connected). This is clear from \S\ref{sss:elliptic curves} (take $G=A[p^n]$, where $A$ is a suitable elliptic curve over $R=\BF_p[t]$).

\section{The algebraic stack $\BT_n^h$}
\subsection{Definition and algebraicity of $\BT_n^h$}
Recall that a \emph{groupoid} is a category in which all morphisms are invertible.

Fix an integer $h\ge 0$. Given a ring $R$, let $\BT_n^h (R)$ denote the groupoid whose objects are $\BT_n$ groups of height $h$ over $R$ and whose morphisms are isomorphisms of $\BT_n$ groups. The assignment $R\mapsto \BT_n^h (R)$ is a functor. As explained by T.~Wedhorn in \cite[Prop.~1.8]{Wedhorn}, this functor is an algebraic stack of finite type over $\BZ$. In other words, there exists a presentation
\begin{equation} \label{e:straightforward presentation}
\BT_n^h =X/\Gamma ,
\end{equation}
where $X$ is a scheme of finite type over $\BZ$, $\Gamma$ is a smooth groupoid acting on $X$, and $X/\Gamma$ is the quotient stack.
% (the precise meaning of the words ``quotient stack'' is explained in foundational texts on algebraic stacks such as \cite{?,?}). 
Moreover, the proof of algebraicity of $\BT_n^h$ given in 
\cite{Wedhorn} yields a concrete presentation 
\begin{equation} \label{e:Wedhorn's presentation}
\BT_n^h =Y/GL(N),
\end{equation}
where $N:=p^{nh}$ and $Y$ is a locally closed subscheme of an affine space over $\BZ$ of dimension $N^2(N+1)$.

Wedhorn's algebraicity result is useful, but his concrete presentation \eqref{e:Wedhorn's presentation} is not really illuminating, and the number $N^2(N+1)$ is typically very big. The problem of finding a better presentation will be discussed in \S\ref{s:the problem}-\ref{s:Katz} below.

For completeness, let us briefly explain the construction of \eqref{e:Wedhorn's presentation} (which is very straightforward).
Recall that a $\BT_n$ group is a commutative finite group scheme with certain extra properties. A commutative finite group $G$ of order $N=p^{nh}$ is the same as a commutative cocommutative Hopf algebra $A$ of dimension $N$ (as an algebra, $A$ is the coordinate ring of~$G$, and the coproduct in $A$ comes from the group operation in $G$). Let $Y'$ be the scheme parametrizing isomorphism classes of commutative cocommutative Hopf algebras of dimension $N$ with a fixed basis. The group $GL(N)$ acts on $Y'$ by changing the basis. The scheme $Y$ from \eqref{e:Wedhorn's presentation} is a locally closed subscheme of $Y'$.

%\subsection{How one gets \eqref{e:Wedhorn's presentation}}

\subsection{Smoothness of $\BT_n^h$}
According to a deep theorem of Grothendieck, the stack $\BT_n^h$ is smooth over $\BZ$, see \cite{Il}.

\section{The problem}   \label{s:the problem} 
\subsection{The stacks $\BT_{n,m}^{h,d}$} 
Recall that $\BT_n^h$ is an algebraic stack over $\BZ$. Let $\BT_{n,m}^h$ denote the reduction of $\BT_n^h$ modulo~$p^m$. It is well known that
\[
\BT_{n,m}^h=\bigsqcup_{d=0}^{h}\BT_{n,m}^{h,d}\, ,
\]
where $\BT_{n,m}^{h,d}$ is defined as follows: if $R$ is a field of characteristic $p$ then $\BT_{n,m}^{h,d}(R)$ is the groupoid of those $\BT_n$ groups over $R$ whose Lie algebra has dimension $d$.

\subsection{Formulation of the problem}   \label{ss:Formulation of the problem}  

%Technically, $\BT_{n,m}^h (R)=\BT_n}^h(R)$ if $p^m$ equals $0$ in $R$, otherwise $\BT_{n,m}^h (R)=\emptyset$.

The following problem is open but not hopeless.

\begin{problem}   \label{the problem} 
Find an illuminating and  ``Shimurizable'' description of the stack $\BT_{n,m}^{h,d}$  for any~$h,n,m$. 
\end{problem}

Here ``Shimurizable'' means ``admitting a generalization related to Shimura varieties'' (including those related to exceptional groups $E_6,E_7$), see \S\ref{ss:Shimura}
below (especially \S\ref{sss:Expectation}) for details.

\subsection{On Shimura varieties}   \label{ss:Shimura}
\subsubsection{Example}
For $d\in\BN$, let $\cA_d$ denote the stack of principally polarized abelian varieties of dimension $d$. This is a stack over $\BZ$. The stack $\cA_d$ is an example of a Shimura variety 
($\cA_d$ is not far from being a scheme, so the word ``variety'' is not misleading).  The complex analytic variety $\cA_{d,\BC}$ corresponding to $\cA_d$ (which is strictly speaking, a stack) has a nice description via uniformization, namely
\begin{equation}   \label{e:A_d}
\cA_{d,\BC}=H_d/Sp(2d,\BZ),
\end{equation}
where $H_d$ is the Siegel upper half-plane (in particular, $H_1\subset\BC$ is the usual half-plane $\im z>0$). As a real variety, $H_d$ has the following description:
\begin{equation}   \label{e:Siegel}
H_d=K\backslash Sp(2d,\BR),
\end{equation}
where $K$ is the maximal compact subgroup of $Sp(2d,\BR)$ (e.g., in the case $d=1$ one has $K=SO(2,\BR)$).

As explained in \S\ref{sss:abelian schemes}, an abelian scheme gives rise to a $\BT_n$ group. Thus we get a morphsim $\cA_d\to\BT_n^{2d}$ and moreover, a morphism
\begin{equation}   \label{e:A_d to BT}
\cA_d\to\BT_n^{2d,Sp},
\end{equation}
where $\BT_n^{2d,Sp}$ is the stack of $\BT_n$ groups $G$ equipped with a symplectic pairing $G\times G\to\BG_m$ (the pairing comes from the polarization of the abelian scheme).
For any $m\in\BN$, one can reduce \eqref{e:A_d to BT} modulo $p^m$ and get a morphism
\begin{equation}   \label{e:reduction of A_d to BT}
\cA_d\otimes\BZ/p^m\to\BT_{n,m}^{2d,Sp}:=\BT_n^{2d,Sp}\otimes\BZ/p^m.
\end{equation}

\subsubsection{General Shimura varieties}
The description of $\cA_{d,\BC}$ given by \eqref{e:A_d}-\eqref{e:Siegel} is based on the group $Sp (2d)$. One can construct analogs of $\cA_{d,\BC}$ by replacing 
$Sp (2d)$ with semisimple groups over $\BQ$ from a certain class. E.g., this class contains some exceptional groups of types $E_6$ and $E_7$.

Constructing analogs of $\cA_d\otimes\BQ$ is much harder, and constructing analogs of $\cA_d$ is very hard (as far as I know, for some types of semisimple groups, the problem is still open). These constructions are rather indirect. The main problem is that for a general Shimura variety we \emph{do not know a modular interpretation} (i.e., an interpretation as a solution to a moduli problem). E.g., we do not know an $E_7$-analog of the notion of abelian variety.

Thus general Shimura varieties are \emph{much more mysterious} than the stacks $\cA_d$. Nevertheless, I believe in the following ``principle of Shimurian democracy'': \emph{all Shimura varieties are equally good} (even if some of them are less accessible than others).

\subsubsection{Expectation}   \label{sss:Expectation}
Based on the above principle, it is natural to expect an analog of \eqref{e:reduction of A_d to BT} (and maybe of \eqref{e:A_d to BT}) for any Shimura variety with good reduction at $p$.

\section{Status of the problem}
%\subsection{Status of the problem}
Problem~\ref{the problem} has a long history, which I try to explain in \S\ref{ss:history} below. In particular, in the 1950's J.~Dieudonn\'e gave a nice description of
$\BT_n^h (R)$ in the case that $R$ is a perfect field of characteristic $p$ (for $n=1$ it is recalled in \S\ref{sss:L,F,V} below).

Let me now say a few words about my e-print \cite{Dr}. Let $G$ be a reductive group over $\BZ_p$ and $\mu:\BG_m\to G$ a \emph{minuscule} homomorphism (i.e., all weights of the corresponding action of $\BG_m$ on the Lie algebra of $G$ are in $\{0,1,-1\}$). For each such pair $(G,\mu )$ and each $n,m\in\BN$, a certain stack
$\BT_ {n,m} ^{G,\mu}$ is defined in Appendix~C of \cite{Dr}; it depends only on the conjugacy class of~$\mu$. It is conjectured in \cite{Dr} that $\BT_ {n,m} ^{G,\mu}$ is algebraic and that
%A certain stack $\BT_ {n,m} ^G$ is defined in \cite{Dr} for any reductive group $G$ from a certain class (essentially, the same one as in the theory of Shimura varieties; in particular, $E_6$ and $E_7$ are allowed). It is conjectured in \cite{Dr} that $\BT_ {n,m} ^G$ is algebraic and that
%\[
%\BT_ {n,m} ^{GL(h)}=\BT_ {n,m} ^h .
%\]
a suitable version of contravariant Dieudonn\'e theory yields an isomorphism
\[
\BT_ {n,m} ^{h,d}\iso\BT_ {n,m} ^{GL(h),{\mu_d}},
%\BT_ {n,m} ^h =\bigsqcup_{\mu\in S}  \BT_ {n,m} ^{GL(h),\mu}.
\]
%where $S$ is the set of conjugacy classes of minuscule homomorp\-hisms~$\BG_m\to PGL(h)$. 
where $\mu_d:\BG_m\to GL(h)$ takes $t\in\BG_m$ to the diagonal matrix whose first $h-d$ diagonal entries equal $t$ and the rest equal $1$. Moreover, in the case $n=m=1$ algebraicity of $\BT_ {n,m} ^{G,\mu}$ is proved in \cite{Dr}.

Note that pairs $(G,\mu )$ with $\mu$ minuscule appear in the theory of Shimura varieties. There exists such a pair in which $G$ has type $E_7$ and $\mu$ is nontrivial.

There are reasons to believe that the conjectures from \cite{Dr} are true and ``provable''. 
The main source of hope is the theory of \emph{prismatic cohomology} developed recently by B.~Bhatt jointly with P.~Scholze and translated into the language of stacks by B.~Bhatt and J.~Lurie. (For this theory, see Bhatt's lecture notes \cite{B}, the references therein, and my lecture \cite{My IAS talk}.)

In fact, the prismatic theory is used already in the definition of $\BT_ {n,m} ^{G,\mu}$. However, as explained in \cite{Dr}, if $R$ is smooth over $\BF_p$ then $\BT_{1,1}^{G,\mu}(R)$ has an elementary description, which goes back to the paper \cite{K72} by N.~Katz. We give it in the next section (for simplicity, in the case $G=GL(h)$).
%in  \S\ref{s:Katz}

Let us note that the stack $\BT_ {n,m} ^{G,\mu}$ is conjectured to be smooth over $\BZ/p^m$. If so then $\BT_ {n,m} ^{G,\mu}$ is uniquely determined by the restriction of the functor $R\mapsto \BT_ {n,m} ^{G,\mu}(R)$ to the category of smooth algebras over $\BZ/p^m$.

%Should I use a different notation for $\BT_ {n,m} ^G$ and explain its relation with the $\BT_ {n,m} ^G$ from \cite{Dr}?

J.~Ansch\"utz and A.~Le Bras have already constructed in \cite{ALB} an isomorphism
\[
\BT_ {\infty,m} ^{h,d}(R)\iso\BT_ {\infty,m} ^{GL(h),{\mu_d}}(R)
\]
assuming that $R$ is quasi-syntomic (here $\BT_ {\infty,m}$ is the projective limit of $\BT_ {n,m}$).

\section{Elementary description of $\BT_{1,1}^{GL(h),\mu_d}$}   \label{s:Katz}
%\subsection{Minuscule homomorphisms $\mu :\BG_m\to PGL(h)$}
\subsection{Points over perfect fields}
\subsubsection{}   \label{sss:L,F,V}
Let $k$ be a perfect field of characteristic $p$. Then $\BT_{1,1}^{GL(h),\mu_d}(k)$ is the groupoid of triples $(L,F,V)$, where $L$ is a vector space over $k$, $F:L\to L$ is a $p$-linear operator, $V:L\to L$ is a $p^{-1}$-linear operator and
\[
\Ker F=\im V, \quad \Ker V=\im F, \quad \dim L=h, \quad \dim\Ker F=d.
\]
The equivalence $\BT_ {1,1} ^{h,d}(k)\iso\BT_ {1,1} ^{GL(h),{\mu_d}}(k)$ is provided by the classical Dieudonn\'e theory.

\subsubsection{}   \label{sss:L,L',L''}

In the situation of \S\ref{sss:L,F,V}, let $L':=\Ker F$, $L'':=\im F$. One can think of an object of $\BT_{1,1}^{GL(h),\mu_d}(k)$ as the following data:

(i) a triple $(L,L',L'')$, where $L$ is a $k$-vector space, $L', L''\subset L$ are subspaces, $\dim L'=d$, $\dim L''=h-d$;

(ii) $p$-linear isomorphisms $L/L'\iso L''$, $L'\iso L/L''$.

\subsection{Points over more general rings}
We are going to give an elementary description of $\BT_{1,1}^{GL(h),\mu_d}(R)$ assuming that $R$ is a smooth algebra over a perfect field of characteristic $p$. 
%(Without any assumptions on $R$, one can give a not so elementary description of $\BT_{1,1}^{GL(h),\mu_d}(R)$ using the magic words ``left Kan extension'', see \cite[Appendix~A]{Dr}).

\subsubsection{}  
Let $R$ be as above.
%Now let $R$ be a smooth $\BF_p$-algebra. 
Let $S:=\Spec R$ and $\Fr:S\to S$ the absolute Frobenius endomorphism.
Then an object of $\BT_{1,1}^{GL(h),\mu_d}(R)$ is the following collection of data:

(i) a triple $(L,L',L'')$, where $L$ is a vector bundle on $S$ and $L', L''\subset L$ are subbundles, $\rank L'=d$, $\rank L''=h-d$;

(ii) isomorphisms 
\begin{equation}  \label{e:horizontal}
\Fr^*(L/L')\iso L'', \quad \Fr^*L'\iso L/L'';
\end{equation}

(iii) an integrable connection $\nabla :L\to L\otimes\Omega^1_S$ such that $\nabla (L'')\subset L''\otimes\Omega^1_S$, the isomorphisms \eqref{e:horizontal} are horizontal (assuming that the connection on the $\Fr$-pullbacks is introduced in the usual way), and the following \emph{Katz condition} holds:
 \begin{equation}   \label{e:Katz condition}
\pCurv_\nabla=-\KS_\nabla ,
 \end{equation}
 where $\pCurv_\nabla$ is the $p$-curvature of $\nabla$ and $\KS_\nabla :L'\to L/L'\otimes\Omega^1_S$  is the map induced by $\nabla :L\to L\otimes\Omega^1_S$.
 The notation $\KS_\nabla$ stands for ``Kodaira-Spencer''; it is motivated by the next remark.
 
\subsubsection{Remark}
Let $\pi :X\to S$ be a smooth proper morphism whose fibers are curves. 
Let $L:=R^1(\pi_*)_{\dR}\cO_X$, where $R^1(\pi_*)_{\dR}$ stands for the first cohomology sheaf of the de Rham direct image.
Then $L$ is equipped with the Gauss-Manin connection $\nabla$, the Hodge filtration $0\subset L'\subset L$, and the conjugate filtration $0\subset L''\subset L$. Moreover, one has horizontal isomorphisms \eqref{e:horizontal}. 
%The surprising equality 
Surprisingly, \eqref{e:Katz condition} holds: this is a particular case of Theorem~3.2 of~\cite{K72}.

\subsubsection{On the Katz condition}
Let us explain why \eqref{e:Katz condition} makes sense.  By definition, $\KS_\nabla$ is a section of $\cA\otimes\Omega^1_S$, where
$\cA:=\HHom(L',L/L')$. Let us show that $\pCurv_\nabla$ is also a section of $\cA\otimes\Omega^1_S$.
\emph{A priori,}  $\pCurv_\nabla$ is a section of $(\Fr_*\EEnd L )^\nabla\otimes\Omega^1_S$, where
$(\Fr_*\EEnd L )^\nabla$ is the sheaf of horizontal sections of $\Fr_*\EEnd L $. 
But the connections on $L'$ and $L/L'$ induced by $\nabla$ are $p$-integrable by virtue of \eqref{e:horizontal}, so $\pCurv_\nabla$ is a section of 
$(\Fr_*\cB)^\nabla\otimes\Omega^1_S$, where $\cB:=\HHom( L / L'', L'')\subset\EEnd L $. 
Finally, $\cB=\Fr^*\cA$ by \eqref{e:horizontal}, so $(\Fr_*\cB)^\nabla=\cA$.

\section{Some history} \label{ss:history}
Problem~\ref{the problem} has a long and complicated history. (However, the current formulation of the problem is recent because the language of algebraic stacks became standard rather recently). I did my best to be as accurate and complete as possible, but I doubt that I succeeded.

The theory developed by J.~Dieudonn\'e in the 1950's yields a description of $\BT_n^h(R)$ in the case that $R$ is a perfect field of characteristic $p$. There is also a similar description if $R$ is any perfect $\BF_p$-algebra, see Theorem~D of \cite{Lau13} and the paragraph after its formulation. Dieudonn\'e's theory is discussed in \cite{Ma, DG, Dem}.

To treat more general rings, Grothendieck suggested in \cite{G1,G2} to develop a \emph{crystalline Dieudonn\'e theory.} Such a theory was developed later in \cite{Me72,MM} and then in the trilogy \cite{BM1,BBM,BM3}. Note that \S 4.4 of \cite{BM3} contains a discussion of problems which remained open in 1990.

Grothendieck proved smoothness of $\BT_n^h$ over $\BZ$. The proof is given in L.~Illusie's paper~\cite{Il}; it is based on Grothendieck's unpublished lectures.

In the 1990's A.~J.~de Jong wrote a series of articles on Barsotti-Tate groups. They are overviewed in his ICM talk  \cite{dJ}.

Around 2000 Th.~Zink started to develop an approach to Barsotti-Tate groups based on his notion of \emph{display.} An overview of some of his works is given in \cite{Me}.

Displays are also used in the works by E.~Lau. His paper \cite{Lau13} was a major source of inspiration for my e-print \cite{Dr}.

Shimurian analogs of Zink's notion of display were introduced by O.~B\"ultel and G.~Pappas in \cite{BP} and then studied in \cite{Lau21}.

\bibliographystyle{alpha}

 \end{document}